\documentclass[preprint,12pt]{article}
\usepackage[width=16cm]{geometry}
\usepackage[english]{babel}
\usepackage[utf8x]{inputenc}
\usepackage{amsmath}
\usepackage{amsthm}
\usepackage{amssymb}
\usepackage{mathtools}
\usepackage{nccmath}
\usepackage{mathcomp}
\usepackage{textcomp}
\usepackage[pdftex]{graphicx}
\usepackage[textwidth=2.2cm]{todonotes}
\usepackage{paralist}
\usepackage{graphics} 
\usepackage{epsfig} 
\usepackage{graphicx}  
\usepackage{epstopdf}
\usepackage{color}
\usepackage{tikz}
\usepackage{abstract} 

\usepackage[colorlinks=true]{hyperref}
\hypersetup{urlcolor=blue, citecolor=red}

\newtheorem{theorem}{Theorem}[section]

\newtheorem{definition}[theorem]{Definition}
\newtheorem{remark}{Remark}

\newcommand{\R}        	{\mathbb {R}}
\newcommand{\mS}		{\mathcal{S}}

\newcommand{\mC}       	{\mathcal {C}}

\title{Local Uniqueness of The Density From Partial Boundary Data for Isotropic Elastodynamics}
\author{Sombuddha Bhattacharyya\footnote{HKUST Jockey Club Institute for Advanced Study, The Hong Kong University of Science and Technology, Clear Water Bay, Kowloon, Hong Kong; email: arkatifr@gmail.com} }\date{}

\begin{document}

\maketitle

\begin{abstract}
We consider an inverse problem in elastodynamics arising in seismic imaging. We prove locally uniqueness of the density of a non-homogeneous, isotropic elastic body from measurements taken on a part of the boundary. We measure the Dirichlet to Neumann map, only on a part of the boundary, corresponding to the isotropic elasticity equation of a 3-dimensional object. 
In earlier works it has been shown that one can determine the sheer and compressional speeds on a neighborhood of the part of the boundary (accessible part) where the measurements have been taken. 
In this article we show that one can determine the density of the medium as well, on a neighborhood of the accessible part of the boundary.
\end{abstract}


\section{Introduction and statement of the main result}\label{Intro}
Let us consider the elastodynamics equation for an isotropic medium. Let $\Omega \subset \R^3$ be a bounded domain with smooth boundary. Consider the following elasticity operator in $\Omega$ defined as
\begin{equation*}
\begin{gathered}
(Pu)_i = \rho\partial^2_tu_i - (Eu)_i \quad
\mbox{where } (Eu)_i = \sum_{jkl = 1}^{3} \partial_{x_j} c_{ijkl} \partial_{x_l} u_k, \qquad i=1,2,3.
\end{gathered}
\end{equation*}
We consider $c_{ijkl} = \lambda \delta_{ij}\delta_{kl} + \mu \delta_{ik} \delta_{jl} + \mu \delta_{il} \delta_{jk}$ to be the elasticity tensor for the isotropic medium and $0<\rho \in C^{\infty}(\Omega)$ is the density with $u=(u_1,u_2,u_3)$ to be the displacement vector.
The elastodynamics equation, for the non-homogeneous isotropic elastic body $\Omega$, can be given as the following initial value problem:
\begin{equation}\label{BVP_1}
\begin{aligned}
&Pu = 0 \qquad && \mbox{in } \Omega \times (0,T),\\
&u|_{\partial \Omega} = f	&& \mbox{for } t \in (0,T),\\
&u|_{t=0} = 0, \quad (\partial_t u)|_{t=0} = 0 && \mbox{in } \Omega.
\end{aligned}
\end{equation}
The principal symbol of the operator $(-E)$ is given by
\begin{equation*}
p_{(-E)}(x,\xi)u = \frac{\lambda + \mu}{\rho} (\xi\cdot u)\xi + \frac{\mu}{\rho}\lvert \xi \rvert^2 u.
\end{equation*}
We readily notice that taking $u = \xi$ or $\xi^{\perp}$ one can get the eigenvalues of $p_{(-E)}(x,\xi)$ as
\begin{equation*}
\begin{aligned}
c_p = \sqrt{\frac{\lambda + \mu}{\rho}}, && c_s = \sqrt{\frac{\mu}{\rho}},
\end{aligned}
\end{equation*} 
with multiplicity $1$ and $2$ respectively corresponding to the eigenspaces $\{r\xi, r\in \R\}$ and $\{\xi^{\perp} \}$.
The $c_p$ and $c_s$ above are known as the speeds of the $p-wave$ and the $s-wave$ respectively.

Now, let us consider the Dirichlet to Neumann map $\Lambda$ on $\partial \Omega \times [0,T]$ as
\begin{equation}\label{DN_map}
(\Lambda f)_i := \sum_{j=1}^{n} \sigma_{ij}(u) \nu^j,
\end{equation}
where $\sigma_{ij}(u) = \lambda (\nabla \cdot u) \delta_{ij} + \mu (\partial_j u_i + \partial_i u_j )$.
In this article we show that the Dirichlet to Neumann map determines the density of the elastic body $\Omega$ on a neighborhood of the boundary where the measurements have been taken.

In \cite{Rachele2000_1} Rachele showed that one can recover $\partial^{m}_{\nu}\rho$, $\partial^{m}_{\nu}\lambda$ and $\partial^{m}_{\nu}\mu$, for $m=0,1,2,\dots$ on the boundary from the Dirichlet to Neumann map $\Lambda$ defined in \eqref{DN_map}. 
Later in \cite{Rachele2000_2,Rachele2003} it has been proved that one can recover the sheer and the compressional speeds $c_s$, $c_p$ and also the density $\rho$ in the domain from the Dirichlet to Neumann map defined on the full boundary. The above results assume that there are no caustics or conjugate points in the domain.
In \cite{Hansen} authors proved the uniqueness of the lens relations and derive several consequences from the Dirichlet to Neumann map for an isotropic elastodynamics with residual stress, without the assumption of caustics and conjugate points.
Recently in \cite{SUV-17-1} the authors have proved that one can recover the sheer and the compressional speeds $c_s$ and $c_p$ in a part of the domain $\Omega$ from the Dirichlet to Neumann map measured only on a part of boundary, assuming some geometric conditions on $\Omega$. 
The geometric condition they consider is the domain (or a part of it) has a strictly convex foliation with respect to the metric $c_{p/s}^{-2}dx^2$. This condition of having a strictly convex foliation admits domains with caustics and conjugate points. For a detailed description of convexly foliated domains see \cite{Monard2014,SUV-17-1,SUV-14}.

We consider geometric optics type solution for the homogeneous isotropic elastodynamics equation and calculate the phase and the amplitude function. We calculate the asymptotic expansion of the amplitude and deduce the necessary conditions for the terms in the expansion of the amplitude in the local coordinate of the the Riemannian metric $g = c_p^{-2}dx^2$.
Simplifying the conditions on the asymptotic expansion and using the given boundary data we get an integral identity involving a 2-tensor consisting of the wave speeds $c_p$, $c_s$ and the density $\rho$.
Hence, our problem reduces to a question of inverting local geodesic ray transform of a symmetric 2-tensors with respect to the Riemannian metric $g = c_p^{-2}dx^2$. 
Here we use the results of \cite{Vasy2016,SUV-14} to invert the ray transform and recover the density on a neighborhood of the boundary where the Dirichlet to Neumann map is given. For more details on geodesic ray transforms of functions and symmetric tensors, see \cite{Sharafutdinov1994,Monard2014,SUV-17-5}.
%
%

Our result determines the density coefficient $\rho$ on a neighborhood of the part of the boundary where the Dirichlet to Neumann map is given. 
One important application of this kind of inverse problem arises in seismic imaging. For more details in the applications of inverse problems in elastodynamics see \cite{Herglotz,Zoeppritz}. 
To define the part of the domain where we prove the uniqueness of the density, let us consider the following.
Let $(\tilde{\Omega},g)$ be an extension of $(\Omega,g)$. 
Let $p \in \partial \Omega$ and $\partial \Omega$ is strictly convex near $p$. Let $\theta \in C^{\infty}(\tilde{\Omega})$ be the local boundary defining function on a neighborhood of $p$ such that $\theta>0$ in $\Omega$. Similar to \cite{Vasy2016,SUV-14} we consider a function $\tilde{x}$ on $\tilde{\Omega}$ as $\tilde{x}(p) = 0$ and $d\tilde{x}(p) = - d\theta(p)$. Fix $c>0$ and consider the neighborhood $\Omega_p = \{ \tilde{x} \geq -c, \theta \geq 0 \}$ of $p$ in $\bar{\Omega}$. Let us denote $\mS = \partial \Omega \cap \partial \Omega_p$. We assume that $\mS$ is strictly convex with respect to $g$ when viewed from $\Omega_p$. The other boundary of the lens shaped domain $\Omega_p$ (i.e. $d\tilde{x} = -c$) is concave when viewed from inside of $\Omega_p$.

\begin{center}
\includegraphics[width=10cm]{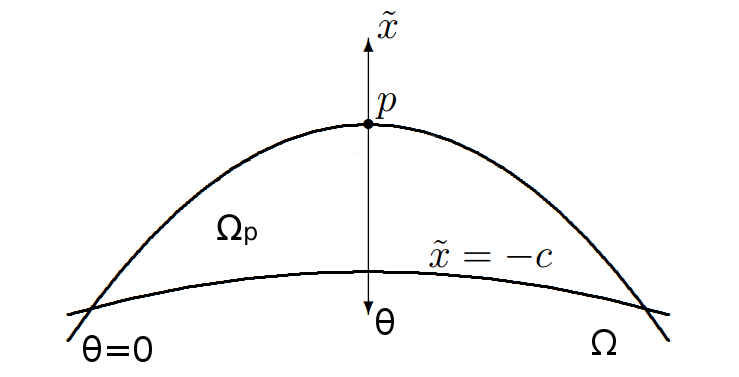}\\
\textbf{Figure 1:} 
\textit{The intersection of $\rho \geq 0$ and $\widetilde{x} > -c$ is the lens shaped region $\Omega_p$.
Note that, as viewed from the inside (i.e. from $\Omega_p$), $\widetilde{x}$ has concave level sets.}
\end{center}
%
%
Having the required notations, let us now state the main theorem of this article.
%

\begin{theorem}\label{Main_Theorem}
Let $\rho_i, \lambda_i, \mu_i$; $i=1,2$ be smooth on $\bar{\Omega}$, satisfying strong convexity condition $3\lambda_i + 2\mu_i > 0$ on $\Omega_p$.
Let $\Lambda_i = \Lambda_2$ on $\mS \times (0,T)$ then $\rho_1 = \rho_2$ on $\Omega_p \setminus \{ x\in \Omega_p : c_p = 2c_s \}$.
\end{theorem}

\begin{remark}
This result gives injectivity of $\rho$ on a neighborhood $\mS$ of $p \in \partial M$ where the Dirichlet to Neumann map is given. The densities are unique on a subset of $\Omega_p$ where $c_p \neq 2c_s$. The assumption of $c_p \neq 2c_s$ also occurs in \cite{Rachele2003} while the author proves the uniqueness of the density from the Dirichlet to Neumann map given on the full boundary.
\end{remark}

\begin{remark}
Under some geometric assumptions on the domain, one can extend the above result in the interior of $\Omega$. One example of this kind of assumption is the strictly convex foliation condition on a neighborhood of $\mS$ in $\Omega$. This condition was used in \cite{SUV-17-1} to prove the uniqueness of the sheer and the compressional speed on a part of $\Omega$ from the Dirichlet to Neumann map given on a part of the boundary.
\end{remark}

In the next section we construct a geometric optics solution of the homogeneous elasticity equation and calculate the asymptotic expansion of the amplitude function involved in the solution. 

\section{Construction of a geometric optics solution}
Consider the following initial boundary value problem
\begin{equation}\label{IVP_1}
\begin{aligned}
&PU = 0,\qquad &&\mbox{in } \Omega \times (0,T),\\
&(U,\partial_t U)|_{t = 0} = (f_0,f_1). \qquad && \mbox{on } \Omega.
\end{aligned}
\end{equation}
%
In this section we construct a geometric optics solution of the above initial value problem \eqref{IVP_1}.
We start with
\begin{equation*}
U = E_0 f_0 + E_1 f_1,
\end{equation*}
where the solution operators $E_k$, $k=0,1$ solves the following equation modulo a smoothing operator.
\begin{equation}\label{E_k PDE}
\begin{aligned}
P E_k &\equiv 0 \quad &&\mbox{on } (0,T) \times \Omega,\\
E_0|_{t=0} &\equiv I, \quad \partial_t E_0|_{t=0} \equiv 0, \quad&& \mbox{on }\Omega,\\
E_1|_{t=0} &\equiv 0, \quad \partial_t E_1|_{t=0} \equiv I, \quad&& \mbox{on }\Omega.
\end{aligned}
\end{equation}
The solution operators, given in terms of Fourier integral operators, are defined as
\begin{equation}\label{E_k FIO}
\begin{aligned}
(E_k v)_j 
=& \int_{\R^3} e^{i\phi^{+}_{p}(t,x,\xi)} a^{j,l}_{+,k,p}(t,x,\xi) \hat{v}_l(\xi)\,d\xi
+ \int_{\R^3} e^{i\phi^{+}_{s}(t,x,\xi)} a^{j,l}_{+,k,s}(t,x,\xi) \hat{v}_l(\xi)\,d\xi\\
&+ \int_{\R^3} e^{i\phi^{-}_{p}(t,x,\xi)} a^{j,l}_{-,k,p}(t,x,\xi) \hat{v}_l(\xi)\,d\xi
+ \int_{\R^3} e^{i\phi^{-}_{s}(t,x,\xi)} a^{j,l}_{-,k,s}(t,x,\xi) \hat{v}_l(\xi)\,d\xi\\
=&\sum_{\pm, p/s, l} \int_{\R^3} e^{i\phi^{\pm}_{p/s}(t,x,\xi)} a^{j,l}_{\pm,k,p/s}(t,x,\xi) \hat{v}_l(\xi)\,d\xi, \quad \mbox{where }v = (v_1,v_1,v_3).
\end{aligned}
\end{equation}
%
The phase function $\phi^{\pm}_{p/s}(t,x,\xi)$ is homogeneous of order 1 in $\xi$ and solves the eikonal equations defined as
\begin{equation}\label{Eqn for phi_1}
\det{p(t,x,\partial_t \phi^{\pm}_{p/s},\nabla_{x}\phi^{\pm}_{p/s})} = 0,
\end{equation}
where $p(t,x,\tau,\xi)$ is the principal symbol of $P$ given as
\begin{equation}\label{symbol p}
p(t,x,\tau,\xi) = -\rho \left[ (\tau^2 - c^2_s \lvert \xi \rvert^2 ) I - (c_p^2- c_s^2)(\xi \otimes \xi) \right].
\end{equation}
One can simplify the condition on $\phi^{\pm}_{p/s}$ given in \eqref{Eqn for phi_1} as
\begin{equation}\label{Eqn for phi_2}
\partial_t \phi^{\pm}_{p/s} = \mp c_{p/s} \lvert \nabla_{x}\phi^{\pm}_{p/s} \rvert.
\end{equation}
We can choose the initial value to be $\phi^{\pm}_{p/s}|_{t=0} = x \cdot \xi$ and solve the above equation using the theory of Hamilton-Jacobi operators.

Unlike the phase function $\phi^{\pm}_{p/s}$, the amplitude cannot be determined using only the principal symbol $p$.
Let us now consider the sub-principal symbol $p_1(t,x,\tau,\xi)$ of $P$ as the sum of the lower order terms in the full symbol of $P$. The expression of $p_1$ is given as
\begin{equation}\label{Sub-principal symbol}
p_1(t,x,\tau,\xi) 
= -i \left[ \nabla_x \lambda \otimes \xi + (\nabla_x \mu \cdot \xi)I + \xi \otimes \nabla_x \mu \right].
\end{equation}
Now, the amplitude $a^{j}_{\pm,k,p/s}(t,x,\xi)$ can be written as
\begin{equation*}
a^{j,l}_{\pm,k,p/s}(t,x,\xi) 
= \sum_{\alpha = 0, -1, -2, \dots}(a^{j,l}_{\pm,k,p/s})_{\alpha},
\end{equation*}
where $(a^{j,l}_{\pm,k,p/s})_{\alpha}$ is homogeneous of order $\alpha$ in $\lvert \xi \rvert$.
In order to calculate the explicit form of $(a^{j,l}_{\pm,k,p/s})_{\alpha}$ let us define the unit vectors $N= \frac{\nabla_x \phi^{\pm}_{p}}{\lvert \nabla_x \phi^{\pm}_{p}\rvert}$ and $N_1$, $N_2$,
are so that $\{N_1, N_2\}$ forms an orthonormal basis of kernel of $p(t,x,\partial_t \phi^{\pm}_{s},\nabla_{x}\phi^{\pm}_{s})$.
Observe that the unit vector $N$ spans the kernel of $p(t,x,\partial_t \phi^{\pm}_{p},\nabla_{x}\phi^{\pm}_{p})$ and $\{N_1, N_2\}$ form an orthonormal basis of kernel of $p(t,x,\partial_t \phi^{\pm}_{s},\nabla_{x}\phi^{\pm}_{s})$.
Now, let us write
\begin{equation}\label{asym_exp}
\begin{aligned}
(a^{j,l}_{\pm,k,p})_{\alpha}
&= (h^{j,l}_{\pm,k,p})_{\alpha} + (b^{l}_{\pm,k,p})_{\alpha}N^{j}, \qquad l =1,2,3,\\
(a^{j,l}_{\pm,k,s})_{\alpha}(t,x,\xi)
&= (h^{j,l}_{\pm,k,s})_{\alpha} 
+ \left[(b^{l}_{\pm,k,s})_{1,\alpha}N^{j}_1 + (b^{l}_{\pm,k,s})_{2,\alpha}N^{j}_2\right],
\end{aligned}
\end{equation}
where $(h^{\cdot,l}_{\pm,k,p/s})_{\alpha}$ is perpendicular to the kernel of $p(t,x,\partial_t \phi^{\pm}_{p/s}, \nabla_x \phi^{\pm}_{p/s})$ for $\alpha \leq -1$ with $(h^{j,l}_{\pm,k,p/s})_{0} = 0$ and $\left(b^{l}_{\pm,k,p}\right)_{\alpha}$, $\left(b^{l}_{\pm,k,s}\right)_{m,\alpha}$, $m=1,2$ are scalars.
Denoting $(a_{\pm,k,p/s})_{\alpha}$ to be the $3 \times 3$ matrix with its $(j,l)$-th coefficients to be  $(a^{j,l}_{\pm,k,p/s})_{\alpha}$ we get
\begin{equation}\label{Eqn for a-alpha}
\begin{aligned}
p(t,x,\partial_t \phi^{\pm}_{p/s},\nabla_{x}\phi^{\pm}_{p/s})(a_{\pm,k,p/s})_{\alpha-1}
= B_{p/s}(a_{\pm,k,p/s})_{\alpha} + C_{p/s}(a_{\pm,k,p/s})_{\alpha+1},
\end{aligned}
\end{equation}
where $(a_{\pm,k,p/s})_1 = 0$ and

\begin{equation}\label{Terms B C}
\begin{aligned}
(B_{p/s}V)
&= i\partial_{\tau,\xi}p\left(t,x,\partial_t \phi^{\pm}_{p/s}, \nabla_x \phi^{\pm}_{p/s}\right) \cdot \partial_{t,x}V
+ i p_1(t,x,\partial_t \phi^{\pm}_{p/s}, \nabla_x \phi^{\pm}_{p/s})V\\
& \qquad \qquad + i\frac{1}{2}\sum_{\lvert \beta \rvert = 2} \partial^{\beta}_{\tau,\xi} p\left(t,x,\partial_t \phi^{\pm}_{p/s}, \nabla_x \phi^{\pm}_{p/s}\right) \cdot \left( \partial^{\beta}_{t,x} \phi^{\pm}_{p/s} \right)V,\\
\left(C_{p/s}V\right)
&= i\partial_{\tau,\xi}p_1\left(t,x,\partial_t \phi^{\pm}_{p/s}, \nabla_x \phi^{\pm}_{p/s}\right) \cdot \partial_{t,x}V\\
& \qquad \qquad + \frac{1}{2}\sum_{\lvert \beta \rvert = 2}\partial^{\beta}_{\tau,\xi}p\left(t,x,\partial_t \phi^{\pm}_{p/s}, \nabla_x \phi^{\pm}_{p/s}\right) \cdot \partial^{\beta}_{t,x} V.
\end{aligned}
\end{equation}
Observe that from \eqref{Eqn for a-alpha} we get a necessary condition that
\begin{equation}\label{Compatibility condition}
N_{p/s} \left[ B_{p/s}(a_{\pm,k,p/s})_{\alpha} + C_{p/s}(a_{\pm,k,p/s})_{\alpha+1} \right] = 0, \quad \mbox{for } \alpha=0,-1,-2,\dots,
\end{equation}
where $N_p = N$ and $N_s = N_1, N_2$.

In the next section we calculate the functions $(a_{\pm,k,p/s})_{\alpha}$ for $\alpha =0,-1$ and get an integral identity relating the Dirichlet to Neumann map.

\section{The equation on the null-bicharacteristics}\label{section_3}
Observe that using \cite{SUV-17-1} one can prove that in $\Omega_p$ the sheer and the compressional speeds ($c_s$ and $c_p$) are unique if the Dirichlet to Neumann map is same on $\mS$. In \cite{SUV-17-1} authors have used the null bicharacteristics which are the curves along which the principal symbol of $P$ is a singular matrix. They proved that the singularities of the boundary data propagates through the null bicharacteristics and from there they proved the result. In this article we go further in to the sub principal symbol of the operator $P$ and show that from the local boundary data we get a local geodesic ray transform involving the density coefficient $\rho$. 

Note that the following calculation is invariant under the choice of different sets of parameters $(\lambda,\mu,\rho)$ hence we remove the index $i$ from the functions $c_p$ and $c_s$ in the following calculation. 
Here we model the path of the null-bicharacteristics of the operator $P$. Observing the equation
\begin{equation*}
0 = \det p(t,x,\tau,\xi) = -\rho^3 (\tau^2 - c_p^2 \lvert \xi \rvert^2)(\tau^2 - c_s^2 \lvert \xi \rvert^2)^2,
\end{equation*}
we get that the null-bicharacteristics of the operator $P$ are along the sets $\{\tau = \pm c_{p/s}\lvert \xi \rvert \}$. Let us consider the case of $p-waves$ in the following analysis.
Consider $\Gamma_{p}^{\pm}$ be the null bicharacteristics corresponding to $\tau = \pm c_p\lvert \xi \rvert$, that is the integral curves of the Hamilton vector fields $V_H$ such that $\det p(t,x,\tau,\xi) = 0$ for $(t,x,\tau,\xi) \in \Gamma_{p}^{\pm}$.
For the Hamiltonian $H = \tau \mp c_p \lvert \xi \rvert$ the vector field $V_H$ is defined by (see \cite{Hormander_1,Rachele2003})
\begin{equation*}
V_H = (\partial_{\tau} H) \partial_t + (\nabla_{\xi} H)\cdot \nabla_x - (\partial_t H) \partial_{\tau} - (\nabla_x H)\cdot \nabla_{\xi}.
\end{equation*}
Hence we get that if $\Gamma_{p}^{\pm}$ is parameterized by $s \in \R$ then $\tau = \pm c_p \lvert \xi \rvert$ and
\begin{equation}\label{path_bichar}
\begin{gathered}
\frac{dt}{ds} = c_p^{-1}, \quad \frac{dx}{ds} = \pm \xi/\lvert \xi \rvert, \quad \frac{d\tau}{ds} = 0, \quad \frac{d\xi}{ds} = \mp \nabla_x(\log c_p) \lvert \xi \rvert,\\
\frac{d}{ds} = c_p^{-1}\partial_t \pm (\xi/\lvert \xi \rvert)\cdot \nabla_x \mp \lvert \xi \rvert \nabla_x(\log c_p)\cdot \nabla_{\xi}.
\end{gathered}
\end{equation}
%
%
Let us consider the compatibility condition \eqref{Compatibility condition} for forward p-wave and $\alpha = 0$. Using equation \eqref{asym_exp} we get
\begin{equation}\label{Eqn_1}
\begin{aligned}
0 &= N B_{p}(a_{\pm,k,p})_{0} = N B_p \left( N \otimes (b_{\pm,k,p})_0\right).
\end{aligned}
\end{equation}
Now to calculate the term $B_p \left(N \otimes (b_{\pm,k,p})_0\right)$ we define the symmetric tensor product as 
\begin{equation*}
u \circledS v = \frac{1}{2}(u \otimes v + v \otimes u).
\end{equation*}
One can show
\begin{equation*}
\begin{aligned}
\partial_{\tau,\xi}p 
&= 2\rho\left(-\tau I,\left[c_s^2\xi_1 I + (c_p^2-c_s^2)(e_1 \circledS \xi)\right], \dots, \left[c_s^2\xi_3 I + (c_p^2-c_s^2)(e_3 \circledS \xi)\right] \right)\\
\partial_{\xi_j\xi_k}p 
&= 2\rho\left[ \delta_{jk}c_s^2I + (c_p^2-c_s^2)(e_j \circledS e_k) \right]\\
N \cdot \partial_{t,x}N &= \frac{1}{2}\partial_{t,x} (N\cdot N) = \frac{1}{2}\partial_{t,x}(1) = 0
\end{aligned}
\end{equation*}
Observe that for compressional waves $N = \frac{\nabla_x \phi^{\pm}_{p}}{\lvert \nabla_x \phi^{\pm}_{p}\rvert}$. 
Writing $(\tau,\xi)=\left(\partial_t \phi^{\pm}_p,\nabla_x \phi^{\pm}_{p}\right)$ and $(b_{\pm,k,p})_0 = b_0$ we get
\begin{equation*}
\begin{aligned}
N \left(\partial_{\tau,\xi}p(t,x,\tau,\xi)\right)\cdot\left(\partial_{t,x}(N\otimes b_0)\right)
&= 2\rho\left[-\tau \partial_t + c_p^2(\xi\cdot\nabla_x) 
+ \frac{1}{2}(c_p^2-c_s^2)\lvert \xi \rvert(\nabla_x \cdot N)
\right]b_0,\\
N p_1(t,x,\tau,\xi)(N\otimes b_0)
&= -i \left[(\xi \cdot \nabla_x)(\lambda + 2 \mu)\right]b_0 = i\rho c_p^2 \left[\xi \cdot \nabla_x (\log \rho c_p^2)\right]b_0. 
\end{aligned}
\end{equation*}
Using \eqref{path_bichar} we get
\begin{equation*}
\begin{aligned}
N \partial_{t}^2 p(t,x,\tau,\xi)\partial_{t}^2(\phi^{\pm}_{p})(N\otimes b_0)
&= 2\rho c_p^2 N(-I)\left[ \pm N(\nabla_x \otimes \xi) N + \frac{d}{ds}(\log c_p) \lvert \xi \rvert \right] (N \otimes b_0)\\
&= -2\rho c_p^2 \left[ \pm N(\nabla_x \otimes \xi) N + \frac{d}{ds}(\log c_p) \lvert \xi \rvert \right] b_0.
\end{aligned}
\end{equation*}
Similarly
\begin{equation*}
\begin{aligned}
N \partial_{\xi_j}\partial_{\xi_k} p(t,x,\tau,\xi)\partial_{x_j}\partial_{x_k}(\phi^{\pm}_{p})(N\otimes b_0)
&= \pm 2\rho N \left[\delta_{jk}c_s^2 I + (c_p^2 - c_s^2) e_j \circledS e_k \right](\partial_{x_j} \xi_k)(N \otimes b_0)\\
&=\pm 2\rho \left[c_s^2 (\partial_{x_j}\xi_j) + (c_p^2 -c_s^2)N_j(\partial_{x_j}\xi_k) N_k\right] b_0
\end{aligned}
\end{equation*}
Hence, from \eqref{Eqn_1} we calculate
\begin{equation*}
\begin{aligned}
0 =& N B_p (N\otimes b_0) \\
=& 2i\rho \left[ -\tau \partial_t + c_p^2(\xi \cdot \nabla_x) + \frac{1}{2}(c_p^2 -c_s^2)\lvert \xi \rvert (\nabla_x \cdot N) + ic_p^2\xi \cdot \nabla_x(\log \rho c_p^2) \right]b_0\\
& - i\rho \left[ c_p^2 N (\nabla_x \otimes \xi)N + c_p^2 \frac{d}{ds}(\log c_p) \lvert \xi \rvert - c_s^2(\nabla_{x}\cdot \xi) - (c_p^2 -c_s^2)N(\nabla_x \otimes \xi)N \right]b_0
\end{aligned}
\end{equation*}
Using \eqref{path_bichar} the above equality reduces to (see \cite[Equation (35)-(37)]{Rachele2003})
\begin{equation*}
\frac{d b_0}{ds} = -\frac{1}{2}\left[\frac{d}{ds}\log(\rho c_p) + (\nabla_x \cdot N) \right]b_0.
\end{equation*} 
Hence for $s_0 \in \R$ so that $(t(s_0), x(s_0)) \in \partial \Omega$ we get
\begin{equation}\label{Solution of ODE}
b_0(s) = b_0(t(s), x(s), \tau, \xi(s)) = b_0(s_0) \sqrt{\frac{\rho c_p (s_0)}{\rho c_p (s)}} \exp \left(-\frac{1}{2} \int_{s_0}^{s} \nabla_x \cdot N\,d\sigma \right).
\end{equation}

Similarly we derive the compatibility condition for $a_{-1} = (a^{j,l}_{\pm,k,p})_{-1}$ given by
\begin{equation*}
N\left[B_p a_{-1} + C_p a_0 \right] = 0.
\end{equation*}
Using the similar calculation as above and from the above equation we get
\begin{equation}\label{ODE a-1}
\begin{aligned}
\frac{d}{ds} a_{-1} + \frac{1}{2}\left[ \frac{d}{ds}(\log \rho c_p) + (\nabla_x \cdot N) \right]a_{-1} = G,\\
\mbox{where }\qquad G = \frac{1}{2i\rho c_p^2 \lvert \xi \rvert}[N B_p h_1 + NC_p(N \otimes a_0)].
\end{aligned}
\end{equation}

Let $\gamma_{\pm}$ be a $p$-wave geodesic, projection of $\Gamma_{p}^{\pm}$ on $M$, with end points $\gamma_{\pm}(s_0) = x_0$ and $\gamma_{\pm}(s_1) = x_1$, where $x_0, x_1 \in \mS$. By our choice of parameterization $\gamma_{\pm}$ is parameterized by $s$. 
Then from the differential equation \eqref{ODE a-1} we get the transport equation
\begin{equation}\label{transport eqn}
g a_{-1} = \int_{\Gamma_{p}^{\pm}} gG + C,
\end{equation}
where $C$ is a constant and 
\begin{equation}\label{step 1}
g = \frac{\sqrt{\rho c_p(x)}}{a_0(s_0) \sqrt{\rho c_p(s_0)}} \exp\left(\frac{1}{2} \int_{s_0}^{s} \nabla_x \cdot N\,ds \right).
\end{equation}
Recall the definitions of $\Omega_p$ and $g$ given in Section \ref{Intro}.
Observe that the analysis above holds locally near $\Gamma_{p}^{\pm}$ whose projection on $(\Omega,g)$ are geodesics. We choose the class of null-bicharacteristics $\Gamma_{p}^{\pm}$ so that the projection of it lies in $\Omega_p$ with end points on $\mS$. 
Using the fact that the Dirichlet to Neumann map is zero on $\mS$ and
a computation of the terms given in equation \eqref{ODE a-1}, \eqref{transport eqn} and \eqref{step 1} similar to \cite[Section 4.1]{Rachele2003} we conclude that 
\begin{equation}\label{Int_Trans}
\int_{\gamma_{\pm}} N \cdot (B N)\,ds = 0,
\end{equation}
where $\gamma_{\pm}$ is a geodesic in $\Omega_p$ with endpoints on $\mS$.
and 
\begin{equation}\label{B_exp}
\begin{aligned}
B(x) 	
=& \frac{A_1 - A_2}{c_p},\\
=& \kappa \left[ \nabla_x \beta_1 \otimes \nabla_x \beta_1 - \nabla_x \beta_2 \otimes \nabla_x \beta_2 \right] - \left[\alpha - \nabla_x\left( \beta_1 - \beta_2 \right) \cdot V\right]I\\
&+ \nabla_x\left(\beta_1 - \beta_2\right) \circledS (V+2\nabla_xc_p) + 2(c_p^2-c_s^2)(\nabla_x \otimes \nabla_x)\left(\beta_1-\beta_2\right).
\end{aligned}
\end{equation}
Where,
\begin{equation*}
\begin{aligned}
\kappa =& \frac{4c_s^2 (c_p^2-2c_s^2)}{c_p(c_p^2-c_s^2)}, \qquad\qquad
\beta_j= \frac{1}{2}\log \rho_j, \qquad \qquad 
V = 2\nabla_x c_p - \frac{8c_s^2}{c_p(c_p^2-c_s^2)}\nabla_x c_s^2,\\
\alpha =& \frac{c_p^2 - 4c_s^2}{2c_p}\Delta_{g}\log \frac{\rho_1}{\rho_2} - \frac{\omega}{4}\left( (\nabla_x \log \rho_1)^2 - (\nabla_x \log \rho_2)^2 \right), \quad
\omega = \frac{c_p^2-4c_s^2}{c_p} + \frac{4c_s^2/c_p}{c_p^2-c_s^2}.
\end{aligned}
\end{equation*}
Using the Theorem 1.1 from \cite{SUV-14} we get $B = dv$ on $\Omega_p$ where $v$ is a smooth 1-form on $\Omega_p$ vanishing on $\mS$.

\section{Inverting the ray transform}
Our aim is to prove $\rho_1 = \rho_2$ on $\Omega_p$ from the integral identity \eqref{Int_Trans}. We already have that $B=dv$ on $\Omega_p$ where $v$ vanishes on $\mS$.
To proceed further we define the \textit{Saint-Venant operator} which will help us to classify the kernel of the geodesic ray transform \eqref{Int_Trans}.
\begin{definition}[Saint-Venant Operator \cite{Sharafutdinov1994}]
The Saint-Venant operator $W$ takes a symmetric 2-tensor to a symmetric 4-tensor defined as
\begin{equation*}
(WB)_{i_1,i_2,j_1,j_2} := \sigma(i_1i_2) \sigma(j_1j_2)\left[ \frac{\partial^2 B_{i_1i_2}}{\partial x_{j_1} \partial x_{j_2}}
- 2\frac{\partial^2 B_{i_1j_1}}{\partial x_{i_2} \partial x_{j_2}} + \frac{\partial^2 B_{j_1j_2}}{\partial x_{i_1} \partial x_{i_2}} \right],
\end{equation*}
where $\sigma(i_1i_2)$ is the symmetrization operator defined for symmetric tensor fields as
\begin{equation*}
\sigma(i_1i_2)V_{i_1i_2i_3i_4} := \frac{1}{2} \sum_{\pi \in \Pi_2} V_{\pi(i_1)\pi(i_2)i_3i_4}, \quad \mbox{where } \Pi_2 \mbox{ is the permutation group of } \{1,2\}.
\end{equation*}
\end{definition}
We denote $d_{g}v = \sigma(\nabla_{g} v)$ (similarly $d_{e}v = \sigma(\nabla_{e} v)$) for any 1-form $v$, where $\nabla_g$ (resp. $\nabla_e$) is the Riemannian connection corresponding to the metric $g$ (the Euclidean metric $e$) and $\sigma$ denotes the symmetrization of the 2-tensor $\nabla_g v$ (or $\nabla_e v$). 
From \cite[Theorem 2.2.1]{Sharafutdinov1994} we see that $0 \equiv W(d_e v)$ for any 1-form $v$ smoothly defined on $\Omega_p$.
%
%
Recall that we have $B=d_{g}v$ in $\Omega_p$ where $v$ is a smooth 1-form on $M$ and $v|_{\mS} = 0$ and $d_g\cdot$ is with respect to the metric $g=c_p^{-2}dx^2$.
Observe that $B$ involves second order derivatives of $\rho_i$ hence $v$ can have at most first order derivatives of $v$.
Moreover $dv = d_{e} v + R(v)$, where $R(v)$ depends on $v$ and derivatives of $c_p$.
Using the \textit{Saint-Venant Operator} we get $W(d_e v) = 0$ that is $W(d_gv) = W(R(v))$ in $\Omega_p$. As $R(v)$ depends only on $v$ and derivatives of $c_p$, so we get that the coefficients of fourth order derivative of $\rho_i$ in $WB$ is zero. 
Let us denote $T_4(A)$ to be the sum of the terms in $A$, having fourth order derivative of $\rho_i$ in $WB$. 

Now, repeating the same arguments in \cite[Section 4.3]{Rachele2003}, for any constant vectors $X, Y$ we get
\begin{equation*}
\begin{aligned}
&W((\nabla_x \otimes \nabla_x) \beta_i) = W(\nabla_x \circledS \nabla_x \beta_i) = 0,\\
&W(\alpha I)(X,X,Y,Y)\\ 
&\qquad \qquad = \lvert X \rvert^2 (Y \cdot \nabla_x)^2 \alpha - 2(X \cdot Y) (X \cdot \nabla_x)(Y \cdot \nabla_x)\alpha + \lvert Y \rvert^2 (X \cdot \nabla_x)^2 \alpha,\\
&\frac{1}{2}W(\nabla_x \beta \otimes \nabla_x \beta)(X,X,Y,Y) \\
&\qquad \qquad = (X \cdot \nabla_x)(Y \cdot \nabla_x)\beta (X \cdot \nabla_x)(Y \cdot \nabla_x)\beta - (X \cdot \nabla_x)^2\beta (Y \cdot \nabla_x)^2\beta
\end{aligned}
\end{equation*}
Therefore for $X = e_i$ and $Y=e_j$, where $e_k \in \R^3$ with $1$ in the $k-$th entry and $0$ otherwise, we have
\begin{equation}\label{T_four}
\begin{aligned}
0 =& T_4(WB) [e_i, e_i,e_j,e_j]\\ 
=& -\sum_{i,j} T_4\left[\kappa W(\nabla_x\beta_1 \otimes \nabla_x\beta_1)
-\kappa W(\nabla_x\beta_2 \otimes \nabla_x\beta_2) - W(\alpha I) \right][e_i, e_i,e_j,e_j]\\
=& -4 T_4 \Delta^2 \alpha\\
=& -2\frac{c_p^2-4c_s^2}{c_p}\Delta^2( \log \rho_1 - \log \rho_2)\\
&+ \frac{c_p^4 - 5c_p^2c_s^2 + 8c_s^4}{c_p(c_p^2 - c_s^2)} \Delta \left[ \nabla_x (\log \rho_1 + \log \rho_2) \cdot \nabla_x (\log \rho_1 - \log \rho_2) \right].
\end{aligned}
\end{equation}
Given the strong convexity condition $3\lambda + 2 \mu >0$ we have $3c_p^2>4c_s^2$ on $\Omega_p$ and hence $\frac{c_p^4 - 5c_p^2c_s^2 + 8c_s^4}{c_p(c_p^2 - c_s^2)}>0$ on $\Omega_p$.
Let us now define $\beta^{\pm} = \log \rho_1 \pm \log \rho_2$, $\gamma = \frac{(c_p^2-c_s^2)(c_p^2 - 4c_s^2)}{c_p^4 - 5c_p^2c_s^2 + 8c_s^4}$ and observe that $\beta^{-}$ solves the fourth order linear elliptic partial differential equation
\begin{equation}
\gamma \Delta^2 \beta^{-} - \Delta \left(\nabla_x \beta^{+} \cdot \nabla_x \beta^{-}\right) = 0, \qquad \mbox{in }\Omega_p.
\end{equation}
Using the result in \cite{Rachele2000_1} one can extend $\rho_1 = \rho_2$, $\lambda_1=\lambda_2$ and $\mu_1 = \mu_2$ smoothly outside $\mS$. 
Note that $\beta^{-}$ vanishes in infinite order in $\mS$.
Let us consider a small open neighborhood $U_0$ of $\mS$ such that $U_0 \cap \Omega \subset \Omega_p \setminus \mC$ where $\mC:= \{ c_p =2c_s \}$.
Choosing $U_0$ suitably we can assume $c_p \neq 2c_s$ in $U_0$ and hence $\gamma$ has a lower bound in $U_0$.
Then for any $x_0 \in U_0 \setminus \bar{\Omega}$ we get
\begin{equation}
\lvert \Delta^2 \beta^{-}(x_0-x)\rvert \leq C\sum_{\lvert \alpha\rvert \leq \beta} \lvert D^{\alpha} \beta^{-}(x_0-x) \rvert, \qquad \mbox{for all }x \in U_0\setminus\{x_0\}.
\end{equation}
By strong unique continuation result of \cite[Section 3]{Protter} we see that $\beta^{-} = 0$ in $U_0$. 
Now varying $U_0$ and use continuity of $\rho_i$ in each component of $\Omega_p \setminus \mC$ we get that $\rho_1\equiv \rho_2$ in $\Omega_p \setminus \mC = \Omega_p \setminus \{ c_p = 2c_s \}$.


\section*{Acknowledgements}
The author expresses his thanks to Professor Gunther Uhlmann for suggesting this problem and for his valuable comments on this project.


\begin{thebibliography}{10}
	
	\bibitem{Hansen}
	S\"onke Hansen and Gunther Uhlmann.
	\newblock Propagation of polarization in elastodynamics with residual stress
	and travel times.
	\newblock {\em Math. Ann.}, 326(3):563--587, 2003.
	
	\bibitem{Herglotz}
	Gustav Herglotz.
	\newblock \"{U}ber die elastizitaet der erde bei beruecksichtigung ihrer
	variablen dichte.
	\newblock {\em Zeitschr. f\"{u}r Math. Phys.}, 52:275--299, 1905.
	
	\bibitem{Hormander_1}
	Lars H\"{o}rmander.
	\newblock {\em The analysis of linear partial differential operators. {I}},
	volume 256 of {\em Grundlehren der Mathematischen Wissenschaften [Fundamental
		Principles of Mathematical Sciences]}.
	\newblock Springer-Verlag, Berlin, second edition, 1990.
	\newblock Distribution theory and Fourier analysis.
	
	\bibitem{Monard2014}
	Fran\c{c}ois Monard.
	\newblock Numerical implementation of geodesic {X}-ray transforms and their
	inversion.
	\newblock {\em SIAM J. Imaging Sci.}, 7(2):1335--1357, 2014.
	
	\bibitem{Protter}
	Murray~H. Protter.
	\newblock Unique continuation for elliptic equations.
	\newblock {\em Trans. Amer. Math. Soc.}, 95:81--91, 1960.
	
	\bibitem{Rachele2000_1}
	Lizabeth~V. Rachele.
	\newblock Boundary determination for an inverse problem in elastodynamics.
	\newblock {\em Comm. Partial Differential Equations}, 25(11-12):1951--1996,
	2000.
	
	\bibitem{Rachele2000_2}
	Lizabeth~V. Rachele.
	\newblock An inverse problem in elastodynamics: uniqueness of the wave speeds
	in the interior.
	\newblock {\em J. Differential Equations}, 162(2):300--325, 2000.
	
	\bibitem{Rachele2003}
	Lizabeth~V. Rachele.
	\newblock Uniqueness of the density in an inverse problem for isotropic
	elastodynamics.
	\newblock {\em Trans. Amer. Math. Soc.}, 355(12):4781--4806, 2003.
	
	\bibitem{Sharafutdinov1994}
	Vladimir~A. Sharafutdinov.
	\newblock {\em Integral geometry of tensor fields}.
	\newblock Inverse and Ill-posed Problems Series. VSP, Utrecht, 1994.
	
	\bibitem{SUV-17-5}
	Plamen Stefanov, Gunther Uhlmann, and Andras Vasy.
	\newblock Local and global boundary rigidity and the geodesic x-ray transform
	in the normal gauge, 2017.
	
	\bibitem{SUV-17-1}
	Plamen Stefanov, Gunther Uhlmann, and Andras Vasy.
	\newblock Local recovery of the compressional and shear speeds from the
	hyperbolic dn map, 2017.
	
	\bibitem{SUV-14}
	Plamen Stefanov, Gunther Uhlmann, and András Vasy.
	\newblock Inverting the local geodesic x-ray transform on tensors, 2014.
	
	\bibitem{Vasy2016}
	Gunther Uhlmann and Andr\'as Vasy.
	\newblock The inverse problem for the local geodesic ray transform.
	\newblock {\em Invent. Math.}, 205(1):83--120, 2016.
	
	\bibitem{Zoeppritz}
	Ernst Wiechert and Zoeppritz Karl.
	\newblock {\em \"{U}ber Erdbebenwellen.}, volume~4.
	\newblock 1907.
	
\end{thebibliography}
\end{document}